\documentclass{article}

\usepackage{hyperref}
\hypersetup{
    colorlinks=true,
    linkcolor=red,
    urlcolor=blue,
    citecolor=green}
\usepackage{arxiv}
\usepackage[utf8]{inputenc} 
\usepackage[T1]{fontenc}    
\usepackage{hyperref}       
\usepackage{url}            
\usepackage{amsfonts}       

\usepackage{amsmath}
\usepackage{graphicx, color}
\usepackage{amssymb}
\usepackage{amsxtra}
\usepackage{amsgen}
\usepackage{amscd}
\usepackage{rotating}
\usepackage{multirow}
\usepackage{latexsym}
\usepackage{fancyhdr}
\usepackage{fontenc}
\usepackage{tikz}
\usetikzlibrary{positioning}
\def\me#1{\ensuremath{E^{#1}_{\phi}}}
\def\df#1#2{\ensuremath{x^{#1}+y^{#1}+z^{#1}\,=\,#2}}
\def\ig#1#2#3{\ensuremath{#1\,#2\,#3}}
\def\ph#1{\ensuremath{#1_{\phi}}}
\usepackage{tikz}
\usepackage{textcomp}
\usepackage{pgfplots}
\pgfplotsset{width=10cm,compat=1.9}
\usepackage{arxiv}
\usepackage[utf8]{inputenc} 
\usepackage[T1]{fontenc}    
\usepackage{hyperref}       
\usepackage{url}            
\usepackage{amsfonts}       
\usepackage{graphicx}
\newtheorem{thm}{Theorem}[section]
\newtheorem{defn}[thm]{Definition}
\newtheorem{eje}[thm]{Example}

\title{Existential Refinement on the search of integer solutions
for the Diophantine equation $x^{3}+y^{3}+z^{3}=n$}

\author{\href{https://orcid.org/0000-0002-4023-2434}{\includegraphics[scale=0.06]{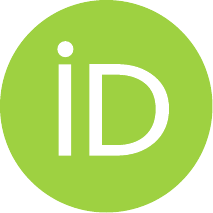}\hspace{1mm}Eduardo Acu\~{n}a} \\
    \texttt{eduardonumberst9@gmail.com} \\
 \AND
    \href{https://orcid.org/0000-0003-4219-2074}{\includegraphics[scale=0.06]{orcid.pdf}\hspace{1mm}Paul Marrero} \\
    \texttt{paulqed0@gmail.com} \\
\AND
\href{https://orcid.org/0000-0001-9026-4883}{\includegraphics[scale=0.06]{orcid.pdf}\hspace{1mm}Samuel
J.Flores}\\
    Department of Mathematics\\
    Carabobo University,FACYT\\
    \texttt{sflores4@uc.edu.ve,\,ssjflores@gmail.com} \\
}

\renewcommand{\shorttitle}{search of integer solutions
for the Diophantine equation $x^{3}+y^{3}+z^{3}=n$}

\begin{document}
\maketitle

\begin{abstract}
We propose a new algorithm, call \textbf{S.A.M} to determinate the
existence of the solutions for the equation \df{3}{n} for a fixed
value \ig{n}{>}{0} unknown.
\end{abstract}


\section{Introduction}
Consider the Diophantine equation
\begin{equation}\label{eq1}
\df{3}{n}
\end{equation}
where $n$ is a fixed positive integer and
$(x,y,z$)$\,\in\,\mathbb{Z}^{3}$, for example: \ig{n}{=}{29} for
inspection the solution is ($x,y,z$)=($3,1,1$); \ig{n}{=}{30} the
solution is discovered in 1999 by E. Pine, K. Yarbrough, W. Tarrant,
and M. Beck \cite{BTPY07}, following an approach suggested by N. Elkies is
($x,y,z$)=($-283059965,-2218888517,2220422932$). Various numerical
investigations (\ref{eq1}) have been carried out, the beginning as
early 1954 \cite{MW55}, for a thorough of the history of these
investigations up to 2000 see \cite{BTPY07}. Heat-Brown \cite{{HB92}} has conjectured
that if $k$ is integer positive with $k\not\equiv\pm4\,(mod\,9)$
then are infinitely triples ($x,y,z$)$\in\mathbb{Z}^{3}$. The paper
of Huisman \cite{{Hui16}} report that solutions are know for a but values of
\ig{k}{<}{100} and $\max\{|x|,|y|,|z|\}\,\leq\,10^{15}$;
\begin{equation}
33,\,42,\,114,\,165,\,390,\,579,\,627,\,633,\,732,\,906,\,921,\,975.
\end{equation}
The computations performed since that the time have been dominated
by an algorithm due to Elkies \cite{{Elk00}}, this algorithm works by finding
rational points near the Fermat curve $X^{3} + Y^{3} = 1$ using
lattice basis reduction; it is well suited to finding solutions for
many values of $k$ simultaneously. But for Hedrick Lecture 1 of
Bjorn Poonen with title Undecidability in number theory \cite{{BP14}}, explain
in the year 2014 that for \ig{n}{=}{33} the result is unknown and
Andrew R. Booker \cite{{AB19}} inspired by the Numberphile video "The
uncracked problem with 33" by Browning and Brady Haran \cite{{BH15}}
\url{https://youtu.be/wymmCdLdPvM} investigate solutions for
\df{3}{k} for few small of $k$ and assume that
$k\equiv\pm3\,(mod\,9)$, $|z|\,>\,\sqrt{k}$ and
$min\{|x|,|y|,|z|\}\,\leq\,10^{16}$ find the the first know
solutions for \ig{k}{=}{33} and \ig{k}{=}{795} found the following:
\begin{eqnarray}
33&=&8866128975287528^{3}+(-8778405442862239)^{3}+(-2736111468807040)^{3}, \nonumber\\
795&=&(-14219049725358227)^{3}+14197965759741571^{3}+2337348783323923^{3}.
\nonumber
\end{eqnarray}
In this paper, in order find existential conditions for calculation
in value fixed $n$ unknown, using the theory of matrix \me{9} \cite{EA}.

\section{OUTLINE NEW SEARCH USING THEORY OF MATRIX \me{9}}
\begin{defn}
Let $v$ a representation decimal of a whole, the digital sum of $v$
\cite{{RL15}} is a function $\phi:\mathbb{Z}\rightarrow\mathbb{Z}$ defined
for:
\begin{equation}
\phi(v)\,=\,
\begin{cases}
0,&\, v\,=\,0, \\
(-1)\sum_{i=0}^{n}|a_{i}|,&\, v\,<\,0, \\
\sum_{i=0}^{n}a_{i},&\, v\,>\,0.
\end{cases}
\end{equation}
\end{defn}

\begin{eje}
\hspace{1cm}
\begin{enumerate}
\item{If $v\,=\,10$, then $\ig{\phi(10)}{=}{\ig{1}{+}{0}}$}\,=\,$1$.
\item{For $v\,=\,-1034$, we have $\phi(-1034)\,=\,(-1)(1+0+3+4)\,=\,-8$.}
\end{enumerate}
\end{eje}


\begin{defn}
Let $v$ be an integer and $\phi$ its digital sum, we say that $\phi$
is a fixed point of v if it happens that:
\begin{equation}
\phi(v)\,=v .
\end{equation}
\end{defn}

\begin{eje}
To everything $v\in [9]_{\geq0}:=\{1,\ldots,9\}$ we have
$\phi(v)\,=\,v.$
\end{eje}

Let $k$ a positive integer and $v$ a positive integer, we denote the
$k-th$ iteration of $\phi$ over $v$ as follows
$\phi^{k}(v)\,=\,\underbrace{\phi\circ\cdots\circ\phi}_{k}(v)$, for
example $\phi^{2}(128)=\phi(\phi(128))=\phi(11)=2$ and $\phi^{0}$ is
the function identity of $\mathbb{Z}$. If $m,\,n$ positive integers,
then $\phi^{m}\circ\phi^{n}\,=\,\phi^{n}\circ\phi^{m}$ and we denote
$v_{m}$ to the integer $v$ that is applied the $m-th$ iteration of
$\phi$, that is, we can represent $v$ as $v_{0}$ and
$\phi^{l}(v)\,=\,v_{l+1}$ with $l\geq0$.
\\
Let us consider the set $\Phi$ formed by all $k\,\geq\,1$ fixed
points of an integer $v$, that is:
\begin{equation}
\Phi_{k}\,=\,\{v\in\mathbb{Z}\,:\,\phi^{k}(v)\,=\,v\}.
\end{equation}
\\
Let $n,\,s$ integers, $[n]_{s}:=[n]\backslash\{s\}$ and
$[n]_{>0}:=\{1,\ldots,n\}$, we denote $\Phi_{1}$ by $D$ where $D$ is
the set of the points fixed of $v\in\,[9]_0$.

\subsection{Condition $k$-interactive of the digital sum}

Let $v$ a representation decimal $d$ and a whole $n$ such what
$\phi$ it's digital sum, the condition $k$-iterative of $v$ is given
as follows:
\begin{enumerate}
\item{If $\phi(v)\,=\,d$ with $d\notin\,[n]_0$, then exists a unique integer $k$
such that $\phi^{k}(v)\,=\,n$ where $n\in\,[9]$}.
\item{If $v\in\,[9]$, then $\phi(v)\,=\,v$}.
\end{enumerate}

\begin{eje}
\hspace{1cm}
\begin{enumerate}
\item{If $\phi(3156)=15,\,\phi(15)\,=\,6$, that is $\phi^{2}(3156)\,=\,6$\,($k=2$).}
\item{If $\phi(4)=4$\,($k=1$).}
\end{enumerate}
\end{eje}


\begin{defn}
Let $n$ be an integer and $k\geq\,1$, the primal class $k$ of $n$,
denoted by $\Phi^{k}$ and is defined as the set formed by the
integers $k$ that are congruent with $n$ module 9, that is,
\begin{equation}
\Phi^{k}\,=\,\{k\in\mathbb{Z}:k\equiv\,n(mod\,9)\}.
\end{equation}
The class primal $k$ of $n$ is on subset of numbers integers is
nonempty where for the axiom of element minimal, exists a unique
element minimal what it we will denote by $Z_{\phi}$.
\end{defn}


\subsection{THEORY OF MATRIX \me{9}}
\begin{defn}
The matrix \cite{EA} \me{9}=($E_{ij,\phi}^{9}$) where $i\,>\,1$ is a
whole, $j\in\,[9]_0$ and it is of the
\begin{eqnarray}
  E_{ij,\phi}^{9} =
  \left\{ \begin{array}{ll}
j, & \textrm{if $i=1$.}\\
9i+j, & \textrm{if $i>1$. }\\
\end{array} \right.
\end{eqnarray}
\end{defn}

\begin{defn}
Let $n\,\in\,\mathbb{Z}$, we defined $k\,\in\,[9]_0$:
\begin{equation}
\phi^{k}\,=\,\{9n+k\}.
\end{equation}
\end{defn}
\label{primaln}
For everything $k\in\,[9]_0$, occurred what $\phi^{k}$, is a
nonempty set of positive integers and assuming the Well Ordering
Principle ($\mathbf{WOP}$)\cite{{NZM91}} , it happens that there exists a minimum element that we will denote by $k_{\phi}$, from so that the columns of the $\me{9}$ are the classes $\phi^{k}$ where
$k\in\,[9]_0$ thus $\me{9}\,=\,(\phi^{k})_{k\in\,[9]_0}$ and what's
more exists a correspondence come in a number integer whole
arbitrary $n$ and its class $\phi^{k}$ defined, choosing a digit
from its fractional part, that is $k_n\,:=\,\{\frac{n}{9}\}$ and
also by construction it happens that $\phi^{k}\,=\,\Phi^{k}$ for
every integer $k$, which implies that the matrix \me{9} is formed by
the primal numbers $Z_{\phi}$ with $Z\in\,[9]_0$.\\

Be $v$ a number whole arbitrary, for decide its class primal
$\phi^{k}$ we calculate $k\,:=\,k_n$ and its primal number that we
denote by $Z_{\phi}$ which is given by $Z_{\phi}\,=\,(k_{n})_{\phi}$
for simplicity we omit the parentheses when we write the primal
number $Z_{\phi}$, for example for $v\,=\,42$, we have that
$42\,=\,6_{\phi}$ and it's class is $\phi^{k_{42}}\,=\,\phi^{6}$ and
to indicate that an integer $v$ is in a class $k$, it is denoted by
$v_{k_{n}}$, if $v\,=\,42$, then we have $42_{6}$.


\subsection{Arithmetic of primal numbers}
The arithmetic used for the refinement and creation of the algorithm comes from the primordial algebra, proved by the "Acuña's Theorem" \cite{EA}.
\subsection{The main theorem}

 The "Acuña's Theorem" \cite{EA}, is the core theorem that expose the principles of the algebra primordial, using the properties of the matrix to arithmetic between classes.
 \subsubsection{notation}
 as mentioned, the theorem that will be exposed below is taken from the Primordial Algebra, being its central theorem, so we will proceed to keep the notations of that article and explain them below:
 \begin{itemize}
     \item $v$ is any integer number
     \item $Z_\phi$ is any Primal number \ref{primaln}.
 \end{itemize}
 When we refer to a $Z_{\phi a}$ and $Z_{\phi a}$ it is any primal number operated with another, the same for the notation of $v_a$, we use the sub-indices $a$ and $b$ to express the notations algebraically.\\
 \vspace{0.2cm}

 \subsection{Acuña's Theorem}
\label{acuña}
 
  For all operation between two or more Primal Classes, the result will be another Primal class, and this is congruent with the numbers inside the class.
 
 \begin{align*}
     \forall \hspace{0.2cm} Z_{\phi a} \odot Z_{\phi b} \equiv  v_a \odot \hspace{0.1cm} v_b~(mod~9)
 \end{align*}
 Where $\odot$ can be any of the following operations: $(+,-,*,/)$\\
 This is a extension of of class of residues by Vinogradov (1977) \cite{Vino}, but this work in both positive and negative integers.\cite{NR}
  
 \subsection{Proof:}
 If:
 \begin{align*}
   \begin{matrix}
    \vspace{0.2cm}
       v_a, Z_{\phi a} \in \phi^a \rightarrow v_a \equiv  Z_{\phi a}~(mod\hspace{0.1cm} 9)\\
       \vspace{0.2cm}
       v_b, Z_{\phi b} \in \phi^b \rightarrow v_b \equiv  Z_{\phi b}~(mod\hspace{0.1cm} \hspace{0.1cm} 9)\\
       v_c, Z_{\phi c} \in \phi^c \rightarrow v_c \equiv  Z_{\phi c}~(mod\hspace{0.1cm} 9)
   \end{matrix}
 \end{align*}
 Then by definition:
 
 \begin{align*}
     \begin{matrix}
      (\phi \odot \phi)^{a+b=}=\left\{ x \in \mathbb Z \mid x \equiv Z_{\phi a} + Z_{\phi b} = Z_{\phi c}~ (mod\hspace{0.1cm} 9) \right\}
     \end{matrix}
 \end{align*}
 From where:
 
 \begin{align*}
     \begin{matrix}
     \vspace{0.2cm}
      x-(Z_{\phi a}+Z_{\phi b}=Z_{\phi c}) = 9t, \hspace{0.2cm} for\hspace{0.1cm} t \in \mathbb Z\\ 
      \vspace{0.2cm}
      x=9t+(Z_{\phi a}+Z_{\phi b}=Z_{\phi c})\\
      \vspace{0.2cm}
      x=9t+[(v_a - 9q)+ (v_b - 9r) = (v_c - 9s)] \hspace{0.2cm}
      for\hspace{0.2cm} q, r, s \in \mathbb Z \\
      \vspace{0.2cm}
      x=9t+[(v_a + v_b)-9(q+r)= v_c - 9s]\\
      \vspace{0.2cm}
      x=9t+[(v_a + v_b) - 9(q+r) + 9s=v_c]\\
      \vspace{0.2cm}
      x=9t+[9(-q-r+s)+(v_a + v_b)=v_c]\\
      \vspace{0.2cm}
      x=9(-q-r+s+t)+(v_a + v_b) = v_c \\
      \rightarrow x \equiv v_a + v_b=v_c \hspace{0.2cm} (mod\hspace{0.1cm} 9) \therefore x \in ( \phi \odot \phi)^{a+b=c}. \hspace{0.2cm} \square
     \end{matrix} 
 \end{align*}
\subsubsection{Addition of primal numbers}
\begin{defn}
Let $\phi^{k}$ where $k\in\,[9]_{\geq\,0}$ primal classes and
$Z_{\phi,\,k}$ it's corresponding number primal of each primal
class, for simplicity we will call $Z_{\phi}$, the addition of
primal numbers consists of several cases that we summarize in the
following table \cite{EA}: \vspace{0,5cm}
\begin{center}
\begin{tabular}{|c|c|c|c|c|c|c|c|c|c|}
  \hline
        + &  $\phi^1$ & $\phi^2$ & $\phi^3$ & $\phi^4$ & $\phi^5$ & $\phi^6$ & $\phi^7$ & $\phi^8$ & $\phi^9$ \\\hline
$\phi^1$ & $2_{\phi}$ & $3_{\phi}$ & $4_{\phi}$ & $5_{\phi}$ &
$6_{\phi}$ & $7_{\phi}$ & $8_{\phi}$ & $9_{\phi}$ &
$1_{\phi}$\\
\hline

$\phi^2$ & $3_{\phi}$ & $4_{\phi}$ & $5_{\phi}$ & $6_{\phi}$ &
$7_{\phi}$ & $8_{\phi}$ & $9_{\phi}$ & $1_{\phi}$ & $2_{\phi}$
\\
\hline

$\phi^3$ & $4_{\phi}$ & $5_{\phi}$ & $6_{\phi}$ & $7_{\phi}$ &
$8_{\phi}$ & $9_{\phi}$ & $1_{\phi}$ & $2_{\phi}$ & $3_{\phi}$
\\
\hline

$\phi^4$ & $5_{\phi}$ & $6_{\phi}$ & $7_{\phi}$ & $8_{\phi}$ &
$9_{\phi}$ & $1_{\phi}$ & $2_{\phi}$ & $3_{\phi}$ & $4_{\phi}$
\\
\hline

$\phi^5$ & $6_{\phi}$ & $7_{\phi}$ & $8_{\phi}$ & $9_{\phi}$ &
$1_{\phi}$ & $2_{\phi}$ & $3_{\phi}$ & $4_{\phi}$ & $5_{\phi}$
\\
\hline

$\phi^6$ & $7_{\phi}$ & $8_{\phi}$ & $9_{\phi}$ & $1_{\phi}$ &
$2_{\phi}$ & $3_{\phi}$ & $4_{\phi}$ & $5_{\phi}$ & $6_{\phi}$
\\
\hline

$\phi^7$ & $8_{\phi}$ & $9_{\phi}$ & $1_{\phi}$ & $2_{\phi}$ &
$3_{\phi}$ & $4_{\phi}$ & $5_{\phi}$ & $6_{\phi}$ & $7_{\phi}$
\\
\hline

$\phi^8$ & $9_{\phi}$ & $1_{\phi}$ & $2_{\phi}$ & $3_{\phi}$ &
$4_{\phi}$ & $5_{\phi}$ & $6_{\phi}$ & $7_{\phi}$ & $8_{\phi}$
\\
\hline

$\phi^9$ & $1_{\phi}$ & $2_{\phi}$ & $3_{\phi}$ & $4_{\phi}$ &
$5_{\phi}$ & $6_{\phi}$ & $7_{\phi}$ & $8_{\phi}$ & $9_{\phi}$
\\\hline
\end{tabular}\\
Table:1 \label{tab1}
\end{center}
\end{defn}

\begin{eje}
\hspace{1cm}
\begin{enumerate}
\item{$18_{9}\,+\,11_{2}\,=\,29_{2}$.}
\item{$12548796_{6}\,+\,67000_{4}\,=\,12615796_{1}$\,($\phi^6\,+\, \phi^4\,=\,1_\phi$).}
\end{enumerate}
\end{eje}

\subsubsection{Subtraction of primal numbers}
\begin{defn}
Let $\phi^{k}$ where $k\in\,[9]_{\geq\,0}$ primal classes and
$Z_{\phi,\,k}$ it's corresponding number primal of each primal
class, for simplicity we will call $Z_{\phi}$, the subtraction of
primal numbers consists of several cases that we summarize in the
following table:
\begin{center}
\begin{tabular}{|c|c|c|c|c|c|c|c|c|c|}
  \hline
        - &  $\phi^1$ & $\phi^2$ & $\phi^3$ & $\phi^4$ & $\phi^5$ & $\phi^6$ & $\phi^7$ & $\phi^8$ & $\phi^9$ \\\hline
$\phi^{-1}$ & $9_{\phi}$ & $1_{\phi}$ & $2_{\phi}$ & $3_{\phi}$ &
$4_{\phi}$ & $5_{\phi}$ & $6_{\phi}$ & $7_{\phi}$ &
$8_{\phi}$\\
\hline

$\phi^{-2}$ & $8_{\phi}$ & $9_{\phi}$ & $1 _{\phi}$ & $2_{\phi}$ &
$3_{\phi}$ & $4_{\phi}$ & $5_{\phi}$ & $6_{\phi}$ & $7_{\phi}$
\\
\hline

$\phi^{-3}$ & $7_{\phi}$ & $8_{\phi}$ & $9_{\phi}$ & $1_{\phi}$ &
$2_{\phi}$ & $3_{\phi}$ & $4_{\phi}$ & $5_{\phi}$ & $6_{\phi}$
\\
\hline

$\phi^{-4}$ & $6_{\phi}$ & $7_{\phi}$ & $8_{\phi}$ & $9_{\phi}$ &
$1_{\phi}$ & $2_{\phi}$ & $3_{\phi}$ & $4_{\phi}$ & $5_{\phi}$
\\
\hline

$\phi^{-5}$ & $5_{\phi}$ & $6_{\phi}$ & $7_{\phi}$ & $8_{\phi}$ &
$9_{\phi}$ & $1_{\phi}$ & $2_{\phi}$ & $3_{\phi}$ & $4_{\phi}$
\\
\hline

$\phi^{-6}$ & $4_{\phi}$ & $5_{\phi}$ & $6_{\phi}$ & $7_{\phi}$ &
$8_{\phi}$ & $9_{\phi}$ & $1_{\phi}$ & $2_{\phi}$ & $3_{\phi}$
\\
\hline

$\phi^{-7}$ & $3_{\phi}$ & $4_{\phi}$ & $5_{\phi}$ & $6_{\phi}$ &
$7_{\phi}$ & $8_{\phi}$ & $9_{\phi}$ & $1_{\phi}$ & $2_{\phi}$
\\
\hline

$\phi^{-8}$ & $2_{\phi}$ & $3_{\phi}$ & $4_{\phi}$ & $5_{\phi}$ &
$6_{\phi}$ & $7_{\phi}$ & $8_{\phi}$ & $9_{\phi}$ & $1_{\phi}$
\\
\hline

$\phi^{-9}$ & $1_{\phi}$ & $2_{\phi}$ & $3_{\phi}$ & $4_{\phi}$ &
$5_{\phi}$ & $6_{\phi}$ & $7_{\phi}$ & $8_{\phi}$ & $9_{\phi}$
\\\hline
\end{tabular}\\
Table: 2 \label{tab2}
\end{center}

\hspace{2cm}
\begin{center}
\begin{tabular}{|c|c|c|c|c|c|c|c|c|c|}
  \hline
        - &  $\phi^{-1}$ & $\phi^{-2}$ & $\phi^{-3}$ & $\phi^{-4}$ & $\phi^{-5}$ & $\phi^{-6}$ & $\phi^{-7}$ & $\phi^{-8}$ & $\phi^{-9}$ \\\hline
$\phi^{1}$ & $-9_{\phi}$ & $-8_{\phi}$ & $-7_{\phi}$ & $-6_{\phi}$ &
$-5_{\phi}$ & $-4_{\phi}$ & $-3_{\phi}$ & $-2_{\phi}$ &
$-1_{\phi}$\\
\hline

$\phi^{2}$ & $-1_{\phi}$ & $-9_{\phi}$ & $-8 _{\phi}$ & $-7_{\phi}$ &
$-6_{\phi}$ & $-5_{\phi}$ & $-4_{\phi}$ & $-3_{\phi}$ & $-2_{\phi}$
\\
\hline

$\phi^{3}$ & $-2_{\phi}$ & $-1_{\phi}$ & $-9_{\phi}$ & $-8_{\phi}$ &
$-7_{\phi}$ & $-6_{\phi}$ & $-5_{\phi}$ & $-4_{\phi}$ & $-3_{\phi}$
\\
\hline

$\phi^{4}$ & $-3_{\phi}$ & $-2_{\phi}$ & $-1_{\phi}$ & $-9_{\phi}$ &
$-8_{\phi}$ & $-7_{\phi}$ & $-6_{\phi}$ & $-5_{\phi}$ & $-4_{\phi}$
\\
\hline

$\phi^{5}$ & $-4_{\phi}$ & $-3_{\phi}$ & $-2_{\phi}$ & $-1_{\phi}$ &
$-9_{\phi}$ & $-8_{\phi}$ & $-7_{\phi}$ & $-6_{\phi}$ & $-5_{\phi}$
\\
\hline

$\phi^{6}$ & $-5_{\phi}$ & $-4_{\phi}$ & $-3_{\phi}$ & $-2_{\phi}$ &
$-1_{\phi}$ & $-9_{\phi}$ & $-8_{\phi}$ & $-7_{\phi}$ & $-6_{\phi}$
\\
\hline

$\phi^{7}$ & $-6_{\phi}$ & $-5_{\phi}$ & $-4_{\phi}$ & $-3_{\phi}$ &
$-2_{\phi}$ & $-1_{\phi}$ & $-9_{\phi}$ & $-8_{\phi}$ & $-7_{\phi}$
\\
\hline

$\phi^{8}$ & $-7_{\phi}$ & $-6_{\phi}$ & $-5_{\phi}$ & $-4_{\phi}$ &
$-3_{\phi}$ & $-2_{\phi}$ & $-1_{\phi}$ & $-9_{\phi}$ & $-8_{\phi}$
\\
\hline

$\phi^{9}$ & $-8_{\phi}$ & $-7_{\phi}$ & $-6_{\phi}$ & $-5_{\phi}$ &
$-4_{\phi}$ & $-3_{\phi}$ & $-2_{\phi}$ & $-1_{\phi}$ & $-9_{\phi}$
\\\hline
\end{tabular}\\
Table: 3 \label{tab3}
\end{center}
\end{defn}

\begin{eje}
\hspace{1cm}
\begin{enumerate}
\item{$200_{2}\,-\,25_{7}\,=\,175_{4}$.}
\item{$20_{2}\,-\,25_{7}\,=\,-5_{-5}$.}
\end{enumerate}
\end{eje}

\subsubsection{Multiplication of primal numbers}

\begin{defn}
Let $\phi^{k}$ where $k\in\,[9]_{\geq\,0}$ primal classes and
$Z_{\phi,\,k}$ it's corresponding number primal of each primal
class, for simplicity we will call $Z_{\phi}$, the multiplication of
primal numbers consists of several cases that we summarize in the
following table:
 
\begin{center}
\begin{tabular}{|c|c|c|c|c|c|c|c|c|c|}
  \hline
        * &  $\phi^1$ & $\phi^2$ & $\phi^3$ & $\phi^4$ & $\phi^5$ & $\phi^6$ & $\phi^7$ & $\phi^8$ & $\phi^9$ \\\hline
$\phi^1$ & $1_{\phi}$ & $2_{\phi}$ & $3_{\phi}$ & $4_{\phi}$ &
$5_{\phi}$ & $6_{\phi}$ & $7_{\phi}$ & $8_{\phi}$ &
$9_{\phi}$\\
\hline

$\phi^2$ & $2_{\phi}$ & $4_{\phi}$ & $6_{\phi}$ & $8_{\phi}$ &
$1_{\phi}$ & $3_{\phi}$ & $5_{\phi}$ & $7_{\phi}$ & $9_{\phi}$
\\
\hline

$\phi^3$ & $3_{\phi}$ & $6_{\phi}$ & $9_{\phi}$ & $3_{\phi}$ &
$6_{\phi}$ & $9_{\phi}$ & $3_{\phi}$ & $6_{\phi}$ & $9_{\phi}$
\\
\hline

$\phi^4$ & $4_{\phi}$ & $8_{\phi}$ & $3_{\phi}$ & $7_{\phi}$ &
$2_{\phi}$ & $6_{\phi}$ & $1_{\phi}$ & $5_{\phi}$ & $9_{\phi}$
\\
\hline

$\phi^5$ & $5_{\phi}$ & $1_{\phi}$ & $6_{\phi}$ & $2_{\phi}$ &
$7_{\phi}$ & $3_{\phi}$ & $8_{\phi}$ & $4_{\phi}$ & $9_{\phi}$
\\
\hline

$\phi^6$ & $6_{\phi}$ & $3_{\phi}$ & $9_{\phi}$ & $6_{\phi}$ &
$3_{\phi}$ & $9_{\phi}$ & $6_{\phi}$ & $3_{\phi}$ & $9_{\phi}$
\\
\hline

$\phi^7$ & $7_{\phi}$ & $5_{\phi}$ & $3_{\phi}$ & $1_{\phi}$ &
$8_{\phi}$ & $6_{\phi}$ & $4_{\phi}$ & $2_{\phi}$ & $9_{\phi}$
\\
\hline

$\phi^8$ & $8_{\phi}$ & $7_{\phi}$ & $6_{\phi}$ & $5_{\phi}$ &
$4_{\phi}$ & $3_{\phi}$ & $2_{\phi}$ & $1_{\phi}$ & $9_{\phi}$
\\
\hline

$\phi^9$ & $9_{\phi}$ & $9_{\phi}$ & $9_{\phi}$ & $9_{\phi}$ &
$9_{\phi}$ & $9_{\phi}$ & $9_{\phi}$ & $9_{\phi}$ & $9_{\phi}$
\\\hline
\end{tabular}\\
Table: 4 \label{tab4}
\end{center}
\end{defn}

\begin{eje}
\hspace{1cm}
\begin{enumerate}
\item{$19_{1}\,*\,31_{4}\,=\,589_{4}$.}
\item{$34_{7}\,*\,81_{9}\,=\,2754_{9}$.}
\end{enumerate}
\end{eje}
\subsection{Division of primal numbers}

\begin{defn}
Let $\phi^{k}$ where $k\in\,[9]_{\geq\,0}$ primal classes and
$Z_{\phi,\,k}$ it's corresponding number primal of each primal
class, for simplicity we will call $Z_{\phi}$, the division of
primal numbers consists of several cases that we summarize in the
following table:
\end{defn}

\begin{center}
\begin{tabular}{|c|c|c|c|c|c|c|c|c|c|}
  \hline
        $\div$ &  $\phi^1$ & $\phi^2$ & $\phi^3$ & $\phi^4$ & $\phi^5$ & $\phi^6$ & $\phi^7$ & $\phi^8$ & $\phi^9$ \\\hline
$\phi^1$ & $1_{\phi}$ & $2_{\phi}$ & $3_{\phi}$ & $4_{\phi}$ &
$5_{\phi}$ & $6_{\phi}$ & $7_{\phi}$ & $8_{\phi}$ &
$9_{\phi}$\\
\hline

$\phi^2$ & $5_{\phi}$ & $1_{\phi}$ & $6_{\phi}$ & $2_{\phi}$ &
$7_{\phi}$ & $3_{\phi}$ & $8_{\phi}$ & $7_{\phi}$ & $9_{\phi}$
\\
\hline

$\phi^3$ & $\varnothing$ & $\varnothing$ & $4_{\phi};1_{\phi};7_{\phi}$ & $\varnothing$ &
$\varnothing$ & $5_{\phi};2_{\phi};8_{\phi}$ & $\varnothing$ & $\varnothing$ & $3_{\phi};9_{\phi};6_{\phi}$
\\
\hline

$\phi^4$ & $7_{\phi}$ & $5_{\phi}$ & $3_{\phi}$ & $1_{\phi}$ &
$8_{\phi}$ & $6_{\phi}$ & $4_{\phi}$ & $2_{\phi}$ & $9_{\phi}$
\\
\hline

$\phi^5$ & $2_{\phi}$ & $4_{\phi}$ & $6_{\phi}$ & $8_{\phi}$ &
$1_{\phi}$ & $3_{\phi}$ & $5_{\phi}$ & $7_{\phi}$ & $9_{\phi}$
\\
\hline

$\phi^6$ & $\varnothing$ & $\varnothing$ & $5_{\phi};2_{\phi};8_{\phi}$ & $\varnothing$ &
$\varnothing$ & $4_{\phi};1_{\phi};7_{\phi}$ & $\varnothing$ & $\varnothing$ & $3_{\phi};9_{\phi};6_{\phi}$
\\
\hline

$\phi^7$ & $4_{\phi}$ & $8_{\phi}$ & $3_{\phi}$ & $7_{\phi}$ &
$2_{\phi}$ & $6_{\phi}$ & $1_{\phi}$ & $5_{\phi}$ & $9_{\phi}$
\\
\hline

$\phi^8$ & $8_{\phi}$ & $7_{\phi}$ & $6_{\phi}$ & $5_{\phi}$ &
$4_{\phi}$ & $3_{\phi}$ & $2_{\phi}$ & $1_{\phi}$ & $9_{\phi}$
\\
\hline

$\phi^9$ & $\varnothing$ & $\varnothing$ & $\varnothing$ & $\varnothing$ &
$\varnothing$ & $\varnothing$ & $\varnothing$ & $\varnothing$ & $Z_{\phi}$
\\\hline
\end{tabular}\\
Table: 5 \label{tab5}
\end{center}

Also we know by the Acuña's Theorem, That:
\begin{align*}
     \frac{\phi^a}{\phi^b}=\phi^c~then~ \frac{\phi^a}{\phi^c}=\phi^b
\end{align*}
Teknomo Kardi, initially had advanced part of the division with digital roots \cite{Kardi}, this table is a complete version.
\subsubsection{Potentiation of primal numbers}
\label{cubo}
\begin{defn}
Let $\phi^{k}$ where $k\in\,[9]_{\geq\,0}$ primal classes and
$Z_{\phi,\,k}$ it's corresponding number primal of each primal
class, for simplicity we will call $Z_{\phi}$, the potentiation of
primal numbers consists of several cases that we summarize in the
following table:
\begin{center}
\begin{tabular}{|c|c|c|c|c|c|c|c|c|c|c|c|c|c|c|}
  \hline
 $()^{*}$ & $()^{2}$ & $()^{3}$ & $()^{4}$ & $()^{5}$ & $()^{6}$ & $()^{7}$ & $()^{8}$ & $()^{9}$ & $()^{10}$ & $()^{11}$ & $()^{12}$ & $()^{13}$ & $()^{14}$ & $()^{15}$ \\
\hline

  $\phi^1$ & $1_{\phi}$ & $1_{\phi}$ & $1_{\phi}$ & $1_{\phi}$ & $1_{\phi}$ & $1_{\phi}$ & $1_{\phi}$ & $1_{\phi}$ & $1_{\phi}$ & $1_{\phi}$ & $1_{\phi}$ & $1_{\phi}$ & $1_{\phi}$ & $1_{\phi}$ \\
\hline

  $\phi^2$ & $4_{\phi}$ & $8_{\phi}$ & $7_{\phi}$ & $5_{\phi}$ & $1_{\phi}$ & $2_{\phi}$ & $4_{\phi}$ & $8_{\phi}$ & $7_{\phi}$ & $5_{\phi}$ & $1_{\phi}$ & $2_{\phi}$ & $4_{\phi}$ & $8_{\phi}$ \\
\hline

  $\phi^3$ & $9_{\phi}$ & $9_{\phi}$ & $9_{\phi}$ & $9_{\phi}$ & $9_{\phi}$ & $9_{\phi}$ & $9_{\phi}$ & $9_{\phi}$ & $9_{\phi}$ & $9_{\phi}$ & $9_{\phi}$ & $9_{\phi}$ & $9_{\phi}$ & $9_{\phi}$ \\
\hline

  $\phi^4$ & $7_{\phi}$ & $1_{\phi}$ & $4_{\phi}$ & $7_{\phi}$ & $1_{\phi}$ & $4_{\phi}$ & $7_{\phi}$ & $1_{\phi}$ & $4_{\phi}$ & $7_{\phi}$ & $1_{\phi}$ & $4_{\phi}$ & $7_{\phi}$ & $1_{\phi}$ \\
\hline

  $\phi^5$ & $7_{\phi}$ & $8_{\phi}$ & $4_{\phi}$ & $2_{\phi}$ & $1_{\phi}$ & $5_{\phi}$ & $7_{\phi}$ & $8_{\phi}$ & $4_{\phi}$ & $2_{\phi}$ & $1_{\phi}$ & $5_{\phi}$ & $7_{\phi}$ & $8_{\phi}$ \\
\hline

  $\phi^6$ & $9_{\phi}$ & $9_{\phi}$ & $9_{\phi}$ & $9_{\phi}$ & $9_{\phi}$ & $9_{\phi}$ & $9_{\phi}$ & $9_{\phi}$ & $9_{\phi}$ & $9_{\phi}$ & $9_{\phi}$ & $9{\phi}$ & $9_{\phi}$ & $9_{\phi}$ \\
\hline

  $\phi^7$ & $4_{\phi}$ & $1_{\phi}$ & $7_{\phi}$ & $4_{\phi}$ & $1_{\phi}$ & $7_{\phi}$ & $4_{\phi}$ & $1_{\phi}$ & $7_{\phi}$ & $4_{\phi}$ & $1_{\phi}$ & $7_{\phi}$ & $4_{\phi}$ & $1_{\phi}$ \\
\hline

  $\phi^8$ & $1_{\phi}$ & $8_{\phi}$ & $1_{\phi}$ & $8{\phi}$ & $1_{\phi}$ & $8_{\phi}$ & $1_{\phi}$ & $8_{\phi}$ & $1_{\phi}$ & $8_{\phi}$ & $1_{\phi}$ & $8_{\phi}$ & $1_{\phi}$ & $8_{\phi}$ \\
\hline

  $\phi^9$ & $9_{\phi}$ & $9_{\phi}$ & $9_{\phi}$ & $9_{\phi}$ & $9_{\phi}$ & $9_{\phi}$ & $9_{\phi}$ & $9_{\phi}$ & $9_{\phi}$ & $9_{\phi}$ & $9_{\phi}$ & $9_{\phi}$ & $9_{\phi}$ & $9_{\phi}$ \\
  \hline
\end{tabular}\\
Table: 6 \label{tab6}
\end{center}
\end{defn}

\begin{eje}
\hspace{1cm}
\begin{enumerate}
\item{($12_{3})^3\,=\,1728_{9}$.}
\item{($45_{9})^{11}\,=\,1532278301220703125_{9}$.}
\item{($2_{2})^8\,=\,256_{4}$.}
\end{enumerate}
\end{eje}

\section{\texorpdfstring{THE ALGORITHM \textbf{S.A.M.} FOR EXISTENCE THE SOLUTIONS OF \df{n}{n}
WITH $n$ FIXED}{}}  
In this section we're interested in determine yourself the existence
of the solutions of The Diophantine equation
$x^{3}+y^{3}+z^{3}\,=\,n$ where $n$ is a positive integer, then we
describe it in the form of an algorithm that supposes a minimum
amount of numerical
restrictions to calculate it:\\

\begin{defn}
Let $\tau:\mathbb{Z}\rightarrow\me{9}$ on function of set, defined
for
\begin{equation}
\tau(n)=\phi^{k},
\end{equation}
where
\begin{equation}
n\equiv\,k\,(mod\,9).
\end{equation}
\end{defn}

\textbf{Input:}  $n$ \\
\textbf{Output:} Exist a solution ($x,\,y,\,z$) of \df{3}{n} or a
message
"nonexistence" if there in no solution. \\
\\
\textbf{step 1:} We calculate the $\tau{(n)}$.\\
\textbf{step 2:} We determine the $\phi^{k}$ defined in the
\textbf{step
1}.\\
\textbf{step 3:} \textbf{If} exist $\phi^{k_{j}}$ where
$k_{j}\in\,[9]_0$ with $j\in\,[3]_{>0}$ that we will denoted by
$X,Y,Z$, that is $\phi^{k_1}=X$,\,$\phi^{k_2}=Y$,\,$\phi^{k_3}=Z$
respectively such that when we calculate $X^{3}\,+\,Y^{3}\,+\,Z^{3}$
and determine their $\phi^{k_3}$ associated such that
$\phi^{k_3}=\phi^{k}$ \textbf{then} existe the solution for
Diophantine equation \textbf{else} output the
message "nonexistence" \textbf{endif}.\\
\\
\textbf{Numerical Example.} If we apply the algorithm \textbf{S.A.M}
for $n\,=\, 33$, we calculate $\tau{(33)}$ which gives us
$\phi^{6}$. On the other hand we determine that
$\phi^{\pm1}=X$,\,$\phi^{\pm8}=Y$,\,$\phi^{\pm9}=Z$, now using the
associated first numbers we have that $X=\pm1_\phi$,\,$Y=\pm8_\phi$
and $Z=\pm9_\phi$, raising to the cube and calculating
$X^{3}\,+\,Y^{3}\,+\,Z^{3}$, using a Punnet combinatorial table, we
determine that the additive combinations of the numbers first, two
by two, we will denote it using letters of the Latin alphabet and
represent them in the following table:

\begin{center}
\begin{tabular}{|c|c|c|c|c|c|c|}
  \hline
  $\oplus$ & $(+1)_{\phi}$ & $(-1)_{\phi}$ & $(+8)_{\phi}$ & $(-8)_{\phi}$ & $(+9)_{\phi}$ & $(-9)_{\phi}$ \\
  \hline
  $(+1)_{\phi}$ & $a=(+2)_{\phi}$   & $b=(+9)_{\phi}$  & $c=(+9)_{\phi}$ & $d=(+2)_{\phi}$ & $e=(+1)_{\phi}$ & $e=(+1)_{\phi}$ \\
  \hline
  $(-1)_{\phi}$ & $b=(+9)_{\phi}$    & $f=(-2)_{\phi}$ & $g=(+7)_{\phi}$ & $h=(-9)_{\phi}$ & $i=(+8)_{\phi}$ & $j=(-1)_{\phi}$ \\
  \hline
  $(+8)_{\phi}$ & $c=(+9)_{\phi}$   & $g=(+7)_{\phi}$ & $g=(+7)_{\phi}$ & $b=(+9)_{\phi}$  & $i=(+8)_{\phi}$ & $i=(+8)_{\phi}$ \\
  \hline
  $(-8)_{\phi}$ & $d=(+2)_{\phi}$   & $k=(-1)_{\phi}$ & $b=(0)_{\phi}$  & $l=(-7)_{\phi}$ & $j=(-1)_{\phi}$ & $i=(+8)_{\phi}$ \\
  \hline
  $(+9)_{\phi}$ & $e=(+1)_{\phi}$   & $i=(+8)_{\phi}$ & $i=(+8)_{\phi}$ & $e=(+1)_{\phi}$ & $d=(+1)_{\phi}$ & $b=(+9)_{\phi}$ \\
  \hline
  $(-9)_{\phi}$ & $e=(+1)_{\phi}$   & $j=(-1)_{\phi}$ & $j=(-1)_{\phi}$ & $m=(-8)_{\phi}$ & $b=(+9)_{\phi}$  & $f=(-2)_{\phi}$ \\
  \hline
\end{tabular}\\
Table: 7 \label{tab7}
\end{center}
\vspace{0,2cm} and the solutions in terms of primal numbers of
Diophantine equation  $x^{3}+y^{3}+z^{3}\,=\,33$, we show them in
the following table:
\begin{center}
\begin{tabular}{|c|c|c|c|c|c|c|c|c|c|c|c|c|c|}
  \hline
   & $a$ & $b$ & $c$ & $d$ & $e$ & $f$ & $g$ & $h$ & $i$ & $j$ & $k$ & $l$ & $m$ \\
  \hline
  $+$ & $\mathbf{\ph{+2}}$ & $\mathbf{\ph{+9}}$ & $\mathbf{\ph{+9}}$ & $\mathbf{\ph{+2}}$ & $\mathbf{\ph{+1}}$ & $\mathbf{\ph{-2}}$ & $\mathbf{\ph{+7}}$ & $\mathbf{\ph{-9}}$ & $\mathbf{\ph{+8}}$ & $\mathbf{\ph{+1}}$ & $\mathbf{\ph{-1}}$ & $\mathbf{\ph{-7}}$ & $\mathbf{\ph{-8}}$ \\
  \hline
  $\mathbf{\ph{+1}}$ & \ph{+3} & \ph{+1} & \ph{+1} & \ph{+3} & \ph{+2} & \ph{+8} & \ph{+8} & \ph{+1} & \ph{+9} & \ph{+2} & \ph{+9} & \ph{+3} & \ph{+2} \\
  \hline
  $\mathbf{\ph{-1}}$ & \ph{+1} & \ph{-1} & \ph{+8} & \ph{+1} & \ph{+9} & \ph{-3} & $\mathbf{\ph{+6}}$ & \ph{-1} & \ph{+7} & \ph{+9} & \ph{-2} & \ph{-8} & \ph{-9} \\
  \hline
  $\mathbf{\ph{+8}}$ & \ph{+1} & \ph{+8} & \ph{+8} & \ph{+1} & \ph{+9} & $\mathbf{\ph{+6}}$ & $\mathbf{\ph{+6}}$ & \ph{+8} & \ph{+7} & \ph{+9} & \ph{+7} & \ph{+1} & \ph{0} \\
  \hline
  $\mathbf{\ph{-8}}$ & \ph{+3} & \ph{-8} & \ph{+1} & \ph{+3} & \ph{+2} & \ph{-1} & \ph{+2} & \ph{-8} & \ph{+9} & \ph{+2} & \ph{-1} & \ph{-6} & \ph{-7} \\
  \hline
  $\mathbf{\ph{+9}}$ & \ph{+2} & \ph{+9} & \ph{+2} & \ph{+2} & \ph{+1} & \ph{+7} & \ph{+7} & \ph{+9} & \ph{+8} & \ph{+1} & \ph{+8} & \ph{+2} & \ph{+1} \\
  \hline
  $\mathbf{\ph{-9}}$ & \ph{+2} & \ph{-9} & \ph{+9} & \ph{+2} & \ph{+8} & \ph{-2} & \ph{+7} & \ph{-2} & \ph{+8} & \ph{+1} & \ph{-1} & \ph{-7} & \ph{-8} \\
  \hline
 \end{tabular}\\
 Table: 8 \label{tab8}
\end{center}
\vspace{1cm} On the other hand, if we calculate the respective
$k_1,\,k_2,\,k_3$ of the triple
\begin{equation}
(x,y,z)=(8866128975287528,-8778405442862239,-2736111468807040)
\end{equation}
described in \cite{{Elk00}} we obtain:
\begin{eqnarray}
k_1\,&=&\,\bigg\{\frac{8866128975287528}{9}\bigg\}\,=\,+2_{\phi}.\\
k_2\,&=&\,\bigg\{\frac{-8778405442862239}{9}\bigg\}\,=\,-7_{\phi}.\\
k_3\,&=&\,\bigg\{\frac{-2736111468807040}{9}\bigg\}\,=\,-4_{\phi}.
\end{eqnarray}

\begin{eqnarray}
(+2_{\phi})^{3}&=&+8_{\phi}.\\
(-7_{\phi})^{3}&=&-1_{\phi}.\\
(-3_{\phi})^{3}&=&-1_{\phi}.\\
\end{eqnarray}

Therefore, from the solutions table, we infer that there are $2$
additional solutions from the one found in \cite{{Elk00}}, totaling three unique combinations for $n\,=\,33$ which
are:
\\
\begin{eqnarray}
\label{19}
(+8_{\phi})\,+\,(-1_{\phi})\,+\,(-1_{\phi})\\
\label{20}
(+8_{\phi})\,+\,(+8_{\phi})\,+\,(+8_{\phi})\\
\label{21}
(-1_{\phi})\,+\,(+8_{\phi})\,+\,(+8_{\phi}).
\end{eqnarray}
\\
The number total of solutions existing (\textbf{s.e}) for $n\,=\,
33$ are $3$.\\
For $n=795$, we determine that $k=3$ and there are $5$ solutions.\\
For $n=1025$, we determine that $k=8$ and there are $8$ solutions.\\
Summarizing we have that:\\
\hspace{1cm}
\begin{center}
\begin{tabular}{c c c}
\hline $n$ &
\textbf{s.e} & \,\,$p$ \\
\hline
33 &   3   &\,\,$\frac{3}{78}$\\
795 &  3   &$\,\,\frac{3}{78}$\\
1025 & 7   &$\,\,\frac{7}{78}$\\
\end{tabular}
\end{center}

where $p$ it means the proportion for the solutions for the $33$,
$795$ and $1025$.

\section{INTERPRETING THE RESULTS OBTAINED}
We can understand by the results of the S.A.M. algorithm the following:\\

\begin{itemize}
    \item There exists a finite set of combinations for each outcome class.
    \item There is a pattern in the proportion of how the outcome classes appear.
\end{itemize}

This shows that there are classes of results with a higher degree of difficulty, while the rest appear more frequently and all this is related to the number of combinations.\\
To calculate these combinations we use what we have learned from the matrix and together with the algorithm manage to obtain a considerable amount of combinations, to optimize the process and have the results fixed, discarding the combinations that are commutative between them and only keep those that are different from each other.\\
In this way it is easier to identify which classes of $k$ will appear more frequently, since they are those that present a greater number of combinations.

\subsection{Combinations for classes k results}
In order to interpret the results we must recapitulate the following:\\
The Acuña's theorem \ref{acuña}, uses the module 9 classification system of the \me{9} and shows us that each operation between such classes will give another congruent class as a result, encompassing both positive and negative classes.\\
This gives rise to a series of algebraic tables that show us the recurrence of appearance of classes referring to each operation, for the one that in this case we care about is that of powers \ref{cubo}, specifically the section of powers cubed in the third column.\\
In this column we can observe that the congruence of cubed classes is {$1_\phi, 8_\phi, 9_\phi$}, as we know that there are negative cubed numbers we obtain the set of congruent classes in this way:
\begin{align*}
    C_\phi:\{-1_\phi,-8_\phi,-9_\phi,~1_\phi,~8_\phi,~9_\phi\}
\end{align*}
Where $C_\phi$ is the subset of classes raised to the cube and these are the elements that the S.A.M. algorithm uses to give us all the possible combinations for $k$. We will call this finite set of solutions the $k_\phi$ classes.\\
The following are the combinations for each congruent $k_\phi$ class:
\subsubsection{\texorpdfstring{Congruent combinations for $k_\phi=1_\phi$}{}}
Replacing X, Y, Z, with the elements of  \{$C_\phi$\} so that their sum results in $1_\phi$ we obtain the following combinations:
\begin{itemize}
    \item $(1_{\phi})+ (1_{\phi})+ (8_{\phi})$
    \item $(1_{\phi})+(1_{\phi})+(-1_{\phi})$
    \item $(-8_{\phi})+(9_{\phi})+(9_{\phi})$
    \item $(1_{\phi})+(9_{\phi})+(-9_{\phi})$
    \item $(1_{\phi})+(8_{\phi})+(-8_{\phi})$\item $(1_{\phi})+(9_{\phi})+(9_{\phi})$
    \item $(1_{\phi})+(-1_{\phi})+(-8_{\phi})$
\end{itemize}
Total of congruent combinations is: 7

\subsubsection{\texorpdfstring{Congruent combinations for $k_\phi=2_\phi$}{}}
Replacing X, Y, Z, with the elements of  \{$C_\phi$\} so that their sum results in $2_\phi$ we obtain the following combinations:
\begin{itemize}
    \item $(9_{\phi})+ (1_{\phi})+ (1_{\phi})$
    \item $(9_{\phi})+(1_{\phi})+(-8_{\phi})$
    \item $(9_{\phi})+(-8_{\phi})+(-8_{\phi})$
    \item $(1_{\phi})+(-8_{\phi})+(-9_{\phi})$
    \item $(-9_{\phi})+(1_{\phi})+(1_{\phi})$
\end{itemize}
Total of congruent combinations is: 5

\subsubsection{\texorpdfstring{Congruent combinations for $k_\phi=3_\phi$}{}}
Replacing X, Y, Z, with the elements of  \{$C_\phi$\} so that their sum results in $3_\phi$ we obtain the following combinations:
\begin{itemize}
    \item $(1_{\phi})+ (1_{\phi})+ (1_{\phi})$
    \item $(1_{\phi})+(-8_{\phi})+(-8_{\phi})$
    \item $(1_{\phi})+(1_{\phi})+(-8_{\phi})$
\end{itemize}
Total of congruent combinations is: 3

\subsubsection{\texorpdfstring{Congruent combinations for $k_\phi=6_\phi$}{}}
\label{camino6}
Replacing X, Y, Z, with the elements of  \{$C_\phi$\} so that their sum results in $6_\phi$ we obtain the following combinations:
\begin{itemize}
    \item $(8_{\phi})+ (-1_{\phi})+ (-1_{\phi})$
    \item $(-1_{\phi})+(8_{\phi})+(8_{\phi})$
    \item $(8_{\phi})+(8_{\phi})+(8_{\phi})$
\end{itemize}
Total of congruent combinations is: 3
\subsubsection{\texorpdfstring{Congruent combinations for $k_\phi=7_\phi$}{}}
Replacing X, Y, Z, with the elements of  \{$C_\phi$\} so that their sum results in $7_\phi$ we obtain the following combinations:
\begin{itemize}
    \item $(9_{\phi})+ (-1_{\phi})+ (-1_{\phi})$
    \item $(-9_{\phi})+(-1_{\phi})+(8_{\phi})$
    \item $(9_{\phi})+(8_{\phi})+(-1_{\phi})$
    \item $(-1_{\phi})+(8_{\phi})+(9_{\phi})$
    \item $(-9_{\phi})+(8_{\phi})+(8_{\phi})$
\end{itemize}
Total of congruent combinations is: 5
\subsubsection{\texorpdfstring{Congruent combinations for $k_\phi=8_\phi$}{}}
Replacing X, Y, Z, with the elements of  \{$C_\phi$\} so that their sum results in $8_\phi$ we obtain the following combinations:
\begin{itemize}
    \item $(9_{\phi})+ (9_{\phi})+ (8_{\phi})$
    \item $(1_{\phi})+(8_{\phi})+(8_{\phi})$
    \item $(8_{\phi})+(-8_{\phi})+(8_{\phi})$
    \item $(-1_{\phi})+(9_{\phi})+(8_{\phi})$
    \item $(8_{\phi})+(9_{\phi})+(-9_{\phi})$
    \item $(8_{\phi})+(1_{\phi})+(-1_{\phi})$
    \item $(8_{\phi})+(-1_{\phi})+(-8_{\phi})$
\end{itemize}
Total of congruent combinations is: 7

\subsubsection{\texorpdfstring{Congruent combinations for $k_\phi=9_\phi$}{}}
Replacing X, Y, Z, with the elements of  \{$C_\phi$\} so that their sum results in $9_\phi$ we obtain the following combinations:
\begin{itemize}
    \item $(9_{\phi})+ (9_{\phi})+ (9_{\phi})$
    \item $(9_{\phi})+ (-9_{\phi})+ (-9_{\phi})$
    \item $(9_{\phi})+(8_{\phi})+(1_{\phi})$
    \item $(9_{\phi})+ (-8_{\phi})+ (-1_{\phi})$
    \item $\underbrace{(8_{\phi})+(-8_{\phi})}_{-9_\phi}+(9_{\phi})$
    \item $\underbrace{(8_{\phi})+(-8_{\phi})}_{9_\phi}+(9_{\phi})$
    \item $\underbrace{(1_{\phi})+(-1_{\phi})}_{-9_\phi}+(9_{\phi})$
    \item $\underbrace{(1_{\phi})+(-1_{\phi})}_{9_\phi}+(9_{\phi})$
    \item $\underbrace{(9_{\phi})+(-9_{\phi})}_{-9_\phi}+(9_{\phi})$
    \item $\underbrace{(9_{\phi})+(-9_{\phi})}_{9_\phi}+(9_{\phi})$
\end{itemize}

Total of congruent combinations is: 10\\
In this case we will expose why the subtraction of two numbers of the same class is $\pm 9$, or 0 in the case of subtracting the same number, remember that S.A.M. makes an iterative addition, first adds the first two terms $X^3+Y^3$ and then the last term $Z^3$ all in its primal form.
For that reason we have remarked where these cases occur where a positive or negative can occur, these cases are explained in the Primal Algebra \cite{EA}.

\subsubsection{\texorpdfstring{Congruent combinations for $k_\phi=4_\phi$ and $k_\phi=5_\phi$}{}}
Just as there are no solutions for $k$ if $k/9$ is of remainder 4 or 5 in module 9.\\
There does not exist a combination of three elements of the set $C_\phi$ which added together give a $k_\phi$= 4 or 5.\\
This shows us that the same analysis of the class system is able to tell which classes are not solutions to an equation, simply by combining the elements to observe which classes appear and which do not.

\section{\texorpdfstring{GRAPH OF CONGRUENCE'S OF THE $K_\phi$ CLASSES FOR THE SUM OF THE 3 CUBES}{}}
\begin{center}

\begin{tikzpicture}
\begin{axis}[
	x tick label style={
		/pgf/number format/1000 sep=},
	ylabel=number of combinations,
	enlargelimits=0.05,
	legend style={at={(0.5,-0.1)},
	anchor=north,legend columns=-1},
	ybar interval=0.7,
]
\addplot 
	coordinates {(1,7) (2,5)
		 (3,3) (4,0) (5,0) (6,3) (7,5) (8,7) (9,10) (10,0)};

\legend{$Classes~\Phi^k$}
\label{figure1}
\end{axis}
\end{tikzpicture}
figure 1
\end{center}

In this graph \ref{figure1} we can observe how the unique combinations that are results of S.A.M. are symmetrically distributed. The uniformity of this generates a valley with its highest peak in the class of $9_\phi$, with the highest number of combinations and its bottom in classes $4_\phi$ and $5_\phi$ that do not have any combination, therefore no number of that class can be a solution.\\

Another characteristic fact is that among the classes that have the same number of combinations, these combinations are parallel to each other, or share similarities that show them as the opposite version.\\

The combinations with the least representatives are those of the $3_\phi$ and $6_\phi$, with only three unique combinations each, this tendency may be the explanation to the difficulty that this family of numbers have in being solutions for $k$, since the previous solved k=\{33, 42,\} both are from the class $6_phi$ and the missing numbers between $k \leq 1,000$ also belong to the classes of the $3_\phi$ and $6_\phi$.\\
The missing numbers to be found according to Booker and Sutherland\cite{Mordell}(2021), are the following and we will order them according to their classes to expose our example.\\

\begin{center}
\begin{tabular}{c c }
\hline $3_\phi$ & $6_\phi$  \\
\hline
390 &   114   \\
732 &  627   \\
921 & 633   \\
975& \\
\end{tabular}
\end{center}

\subsection{examples with diagrams}
In this section we are going to show how the analysis of the existential algorithm S.A.M. works.\\
First we start from the transformation of an integer to its primal form, this through its digital roots\cite{NR}, then by the Acuaña theorem \ref{acuña}, a whole algebraic set is derived between these classes, of which we are only interested in the cubed classes,
then we will take that subset of classes cubed, positive and negative and combine them into a sum by parts, first we will add the first two classes, then we will add the third class, in the end we will arrive at a single final $k$ class.\\
The objective of S.A.M. is to know the paths of the cubed classes $X,~ Y,~Z$, that satisfy $k$, thus optimizing the search by classes of $k$ and not doing brute force searching for all integer.\\
At the end the path can be seen in a simple diagram as follows:

\begin{eje}

\begin{center}

\begin{tikzpicture}[xyz/.style={circle, draw=green!60, fill=green!5, very thick, minimum size=7mm},
k/.style={circle, draw=cyan!60, fill=cyan!5, very thick, minimum size=7mm},
formula/.style={rectangle, draw=red!60, fill=red!5, very thick, minimum size=5mm},
]
\node[xyz]    (x)   {$X_\phi$};
\node[xyz]    (y)   [right=of x]  {$Y_\phi$};
\node[xyz]    (z)   [right=of y]  {$Z_\phi$};
\node[formula] (suma)  [below=of x] {$+$};
\node[xyz]  (x+y)  [below=of suma] {$X_\phi+Y_\phi$};
\node[formula] (suma2) [right=of x+y] {$+$};
\node[k] (k) [below=of suma2] {$K_\phi$};

\draw[->] (x.south) -- (suma.north);
\draw[->] (y.south) .. controls +(down:7mm) and +(right:7mm) .. (suma.east);
\draw[->] (z.south).. controls +(down:7mm) and +(right:7mm) .. (suma2.east);
\draw[->] (suma.south) -- (x+y.north);
\draw[->] (x+y.east) -- (suma2.west);
\draw[->] (suma2.south) -- (k.north);
\end{tikzpicture}
\label{diagrama1} figure 2
\end{center}
\end{eje}
Having this analysis let's proceed to give an example with the numbers to find the 33, if it belongs to one of the paths already given by the algorithm, as well as giving it in the form of a series of steps to try it with any other number that already knows the result.
\subsubsection{Example with 33}
\emph{STEP 1:}with the numbers obtained by Booker. A. \cite{AB19} we proceed to identify their classes:
\begin{align*}
    X=\underbrace{8866128975287528}_{2_\phi}\\
    Y=\underbrace{-8778405442862239}_{-7_\phi}\\
    Z=\underbrace{-2736111468807040}_{-4_\phi}\\
\end{align*}
\emph{STEP 2:}now with only having the classes of the numbers, we can know what will be its class to the power cubed, so this step is to identify the cubic class, the following diagram simplifies more the power table\ref{cubo} by only showing the cubic part.
\begin{center}
\begin{tikzpicture}[n/.style={circle, draw=black, fill=white, very thick, minimum size=7mm},
clase/.style={circle, draw=black, fill=cyan!10, very thick, minimum size=7mm},
k/.style={circle, draw=green!60, fill=green!5, very thick, minimum size=7mm},
formula/.style={rectangle, draw=red!60, fill=red!5, very thick, minimum size=5mm},
cuadro/.style={rectangle, draw=black, fill=white, very thick, minimum size=5mm},
]
\node[cuadro]   (primal)   {$\pm~Z_\phi$};
\node[formula]   (cubo) [right=of primal] {$(Z_\phi)^3$};
\node[n]   (1)    [below=of primal]   {$\pm1_\phi$};
\node[k]   (11)  [right=of 1]  {$\pm~1_\phi^`$};
\node[clase]   (2)    [below=of 1]   {$\pm~2_\phi$};
\node[k]   (81)  [right=of 2]  {$\pm~8_\phi^`$};
\node[n]   (3)    [below=of 2]   {$\pm~3_\phi$};
\node[k]   (91)  [right=of 3]  {$\pm~9_\phi^`$};
\node[clase]   (4)    [below=of 3]   {$\pm~4_\phi$};
\node[k]   (12)  [right=of 4]  {$\pm~1_\phi^`$};
\node[n]   (5)    [below=of 4]   {$\pm~5_\phi$};
\node[k]   (82)  [right=of 5]  {$\pm~8_\phi^`$};
\node[n]   (6)    [below=of 5]   {$\pm~6_\phi$};
\node[k]   (92)  [right=of 6]  {$\pm~9_\phi^`$};
\node[clase]   (7)    [below=of 6]   {$\pm~7_\phi$};
\node[k]   (13)  [right=of 7]  {$\pm~1_\phi^`$};
\node[n]   (8)    [below=of 7]   {$\pm~8_\phi$};
\node[k]   (83)  [right=of 8]  {$\pm~8_\phi^`$};
\node[n]   (9)    [below=of 8]   {$\pm~9_\phi$};
\node[k]   (93)  [right=of 9]  {$\pm~9_\phi^`$};
\draw[->] (primal.east) to node[above] {$(~)^3$} (cubo.west);
\draw[->] (1.east) -- (11.west);
\draw[->] (2.east) -- (81.west);
\draw[->] (3.east) -- (91.west);
\draw[->] (4.east) -- (12.west);
\draw[->] (5.east) -- (82.west);
\draw[->] (6.east) -- (92.west);
\draw[->] (7.east) -- (13.west);
\draw[->] (8.east) -- (83.west);
\draw[->] (9.east) -- (93.west);
    \end{tikzpicture}
figure 3
\end{center}
With the previous diagram we identified that the path for 33 will be as follows: $Z^3=8_\phi,~X^3=-1_\phi,~Y^3=-1_\phi$.\\
With this step we save ourselves from raising to the cube since we know in which class it will end up, but to make the example more precise we will show the exact numbers that added together give 33 and their classes.
\begin{align*}
    X^3=\underbrace{696950821015779435648178972565490929714876221952}_{8_\phi}\\
    Y^3=\underbrace{-676467453392982277424361019810585360331722557919}_{-1_\phi}\\
    Z^3=\underbrace{20483367622797158223817952754905569383153664000}_{-1_\phi}
\end{align*}
\emph{STEP 3:}In this step we simply replace $X$, $Y$, $Z$, by their classes and proceed to operate, if they satisfy the class of the result, then $X_\phi,~ Y_\phi,~ Z_\phi,$ is one of the paths for that class.

\begin{center}
    \begin{tikzpicture}[xyz/.style={circle, draw=green!60, fill=green!5, very thick, minimum size=7mm},
k/.style={circle, draw=cyan!60, fill=cyan!5, very thick, minimum size=7mm},
formula/.style={rectangle, draw=red!60, fill=red!5, very thick, minimum size=5mm},
]
\node[xyz]    (x)   {$-1_\phi$};
\node[xyz]    (y)   [right=of x]  {$-1_\phi$};
\node[xyz]    (z)   [right=of y]  {$8_\phi$};
\node[formula] (suma)  [below=of x] {$+$};
\node[xyz]  (x+y)  [below=of suma] {$-2_\phi$};
\node[formula] (suma2) [right=of x+y] {$+$};
\node[k] (k) [below=of suma2] {$6_\phi$};

\draw[->] (x.south) -- (suma.north);
\draw[->] (y.south) .. controls +(down:7mm) and +(right:7mm) .. (suma.east);
\draw[->] (z.south).. controls +(down:7mm) and +(right:7mm) .. (suma2.east);
\draw[->] (suma.south) -- (x+y.north);
\draw[->] (x+y.east) -- (suma2.west);
\draw[->] (suma2.south) -- (k.north);
\end{tikzpicture}
figure 4
\end{center}
As can be seen, the combination of classes for 33 is one of the three possible paths already found by the algorithm \ref{camino6}, following these steps we analyzed each number found that satisfies $K$ from $K \leq 1000$.\\
Also other numbers that share the same class as $6_\phi$ have this path, you can verify it by following these steps with the following \href{https://arxiv.org/src/1604.07746v1/anc/sumofthreecubes_20160426.txt}{data} found in \cite{Hui16} 
\section{Conclusions}
In terms of optimization, what the S.A.M. algorithm achieves is to be able to determine the paths that exist to find a given class for a $K$.\\
It also identifies that all the numbers that have been the most elusive for the sum of 3 cubes belong to the classes of $3_\phi$ and $6_\phi$, those that as has been shown have fewer combinations, therefore, have a smaller range of occurrence.\\
Since there is more than one combination for each class, it is possible that there is more than one way to calculate the same number.\\
\subsection{thoughts about the process}

The primordial algebra is a new tool and a new mathematics, often during the research, we found that we had to build from scratch all the blocks that compose it, showing that there are gaps in modern mathematics that are missing to be filled, due to this although its application is simple its argumentation is very abstract, this is the reason why we have tried to make as clear as possible the whole process of analysis.
Actually problems like this require different approaches, mathematics is a game where you create your own pieces as long as they are true, in the middle of this logic crisis it can be complicated to see new approaches to which we are used to, we look for the mystery of numbers, in simple things there can be complicated arguments, it is necessary to relearn from another point of view to solve the mysteries that are hidden in our mind.

\section{ACKNOWLEDGEMENTS}
We thank Juan Yaguaro for the discussions with the computational
codes numerical, this work was carried out using the computational
facilities of the Department of Mathematics, Andrew Booker for their comments, FACYT and
'$I\alpha\acute{\omega}$ for your guide.

\bibliographystyle{unsrtnat}


\begin{thebibliography}{1}

    \bibitem{AB19}
    Andrew. R Booker.
    \newblock \textit{CRACKING THE PROBLEM WITH 33},
    \newblock preprint arXiv:1804.09028, 2019.

    \bibitem{BP14}
    Bjorn Poonen.
    \newblock \textit{Undecidability in number theory}.
    \newblock UAM Colloquium March 25, 2010.

    \bibitem{BH15}
    Tim Browning and Brady Haran.
    \newblock \textit{The uncracked problem with 33}.
    \newblock, 2015, \url{https://youtube/wymmCdLdPvM}.

    \bibitem{BTPY07}
     Michael Beck, Erick Pine, Wayne Tarrant and Kim Yarbrough
    \newblock \textit{New integer representations as the sum of three cubes}.
    \newblock
    \textnormal{
    Math.\,Comp.\,\textbf{76}\,(2007),\,no.\,259,\,16831690.\,MR2299795}

    \bibitem{Elk00}
    Noam D. Elkies.
    \newblock \textit{Rational points near curves and small nonzero $|x^3\,-\, y^2|$ via lattice reduction}
    \newblock  \textnormal{Algorithmic number theory (Leiden, 2000), Lecture Notes in Comput.
    Sci.}
    \newblock \textnormal{vol.\,1838,\,Springer,\, Berlin,\, 2000,\, pp.\, 3363.\, MR
    1850598}

    \bibitem{HB92}
    D.R Heath-Brown.
    \newblock \textit{The density of zeros of form for which weak
    approximations fails}.
    \newblock Math. Comp. \textbf{59}\, (1992),\,no.\, 200,\, 613623. MR 1146835

    \bibitem{Hui16}
    Sander G. Huisman.
    \newblock \textit{Newer sums of three cubes}.
    \newblock \textnormal{preprint arXiv:1604.07746, 2016}.

    \bibitem{MW55}
    J. C. P. Miller and M. F. C. Woollett.
    \newblock \textit{Solutions of the Diophantine equation} $x^3+y^3+z^3=k$.
    \newblock \textnormal{J. London Math.\, Soc.\, 30\, (1955),\, 101110.\, MR
    00679}

    \bibitem{NZM91}
    Ivan Niven, Herbert S. Zuckerman, and Hugh L. Montgomery.
    \newblock \textit{An introduction to the theory of numbers}.
    \newblock \textnormal{fifth ed., John Wiley and Sons,\, Inc.,\, New York,\, 1991.\, MR
    1083765}

    \bibitem{RL15}
    Ram\'irez L\'opez, Gorrostola Nadad.
    \newblock \textit{Suma Sucesiva}.
    \newblock \textnormal{2015,\, Universidad del Atl\'antico p.\,52.}
    
    \bibitem{EA}
    Romero, Paul Francisco Marrero. "Matrix$E_\phi^ 9$. Classification of integers and primordial algebra." MATUA. Revista de matemáticas de la Universidad del Atlántico (Colombia). 8.1 (2021): 10-45.

    \bibitem{NR}
    P. F. M. Romero and E. J. A. T. ABOUT THE NEGATIVE DIGITAL ROOT AND SOME OF ITS PROPERTIESRELATED TO MODULAR ARITHMETIC., Oct. 2021. URLhttps://doi.org/10.5281/zenodo.5598339.


    \bibitem{Vino}
    I. M. Vinogradov and E. A. Bernardo.Fundamentos de la teoría de los números. Mir, 1977.
    
    \bibitem{Kardi}
    K. Teknomo, What is Digital Root, \url{https://people.revoledu.com/kardi/}
    
    \bibitem{Mordell}
    Booker, Andrew R., and Andrew V. Sutherland. "On a question of Mordell." Proceedings of the National Academy of Sciences 118.11 (2021).
    
 \end{thebibliography}

\end{document}